\documentclass{article}%
\usepackage{amsmath}
\usepackage{amsfonts}
\usepackage{amssymb}
\usepackage{graphicx}%
\providecommand{\U}[1]{\protect\rule{.1in}{.1in}}
\newtheorem{theorem}{Theorem}

\newenvironment{proof}[1][Proof]{\noindent\textbf{#1.} }{\ \rule{0.5em}{0.5em}}
\begin{document}

\title{Commensurability and bi-interpretability of groups}
\author{Dan Segal}
\maketitle

In \cite{NST}, Section 3 we suggest that the property of being finitely
axiomatizable (in a given class of groups) should be inherited by finite
extensions and by definable subgroups of finite index. The underlying idea is
this: if $H$ is a definable normal subgroup of finite index in a group $G$,
then first-order properties of $H$ can be read as properties of $G$, and
first-order properties of $G$ can be encoded as properties of $H$.

A general framework for expressing such situations is the concept of
\emph{bi-interpretability}, see for examaple \cite{H}, \S \S 5.3, 5.4. The
following observation is probably in the folklore, but seems worth recording:

\begin{theorem}
\label{T1}Let $G$ be a group and $H\neq1$ a definable normal subgroup of
finite index. The following are equivalent:

\begin{description}
\item[a)] for each $t$ in some transversal to $G/H$ the conjugation action of
$t$ on $H$ is definable in $H$.

\item[b)] $G$ is bi-interpretable with $H$.
\end{description}
\end{theorem}

\begin{proof}
Assume \textbf{(a)}. We may suppose that $\left\vert G:H\right\vert =m\geq2$.
Say $H=\kappa(G),$ fix a transversal $\{t\,_{1}=1,t_{2},\ldots,t_{m}\}$ to
$G/H$, and let $\sigma_{i}:H\rightarrow H$ be the definable function
$h\longmapsto t_{i}ht_{i}^{-1}$. Define $c_{ij}\in H$ and $k(i,j)$ by
$t_{i}t_{j}=c_{ij}t_{k(i,j)}$ for each $i$ and $j$.

\emph{Interpreting} $H$ \emph{in} $G$. The inclusion map $H\rightarrow G$
identifies $H$ with $\kappa(G)$. The goup operation is inherited from $G$.

\emph{Interpreting} $G$ \emph{in} $H$. Fix an element $\xi\in H\smallsetminus
\{1\}$. Suppose first that $m\geq3$. Define subsets $\Gamma_{i,0},\Gamma
_{i,1}\subset H^{m}$ as follows:%
\begin{align*}
\gamma &  \in\Gamma_{i,0}\Longleftrightarrow\gamma_{i}\neq1\wedge
\bigwedge_{j\neq i}\gamma_{j}=1\\
\gamma &  \in\Gamma_{i,1}\Longleftrightarrow\gamma_{i}=1\wedge\bigwedge_{j\neq
i}\gamma_{j}=\xi.
\end{align*}
Then each $\Gamma_{i,0}$ and $\Gamma_{i,1}$ is a definable subset of $H^{m}$,
as is $\Gamma:=\bigcup_{i=1}^{m}(\Gamma_{i,0}\cup\Gamma_{i,1})$.

Suppose $g=ht_{i}\in G$. Define $\gamma(g)=(\gamma_{1},\ldots,\gamma_{m})$ by
\begin{align*}
\gamma_{i}  &  =h,~\ \gamma_{j}=1~(j\neq i)\text{ \ if }h\neq1\\
\gamma_{i}  &  =1,~\ \gamma_{j}=\xi~(j\neq i)\text{ \ if }h=1.
\end{align*}
Thus $\gamma(g)\in\Gamma_{i,0}\cup\Gamma_{i,1}$. In fact $g\longmapsto
\gamma(g)$ identifies $(H\smallsetminus\{1\})t_{i}$ with $\Gamma_{i,0}$ and
$\{1.t_{i}\}$ with $\Gamma_{i,1}$.

We define a binary operation on $\Gamma$ so as to make this mapping $\gamma$ a
group isomorphism.

Thus on $\Gamma_{i,0}\times\Gamma_{j,0}$:%
\[
(1,\ldots,h,\ldots,1)\cdot(1,\ldots,k,\ldots,1)=\left\{
\begin{array}
[c]{ccc}%
(1,\ldots,hk^{\sigma_{i}}c_{ij},\ldots,1)\in\Gamma_{k(i,j),0} & \text{if} &
hk^{\sigma_{i}}c_{ij}\neq1\\
&  & \\
(\xi,\ldots,1,\ldots,\xi)\in\Gamma_{k(i,j),1} & \text{if} & hk^{\sigma_{i}%
}c_{ij}=1
\end{array}
\right.  ;
\]

on $\Gamma_{i,1}\times\Gamma_{j,0}$:%
\[
(\xi,\ldots,1,\ldots\xi)\cdot(1,\ldots,k,\ldots,1)=\left\{
\begin{array}
[c]{ccc}%
(1,\ldots,k^{\sigma_{i}}c_{ij},\ldots,1)\in\Gamma_{k(i,j),0} & \text{if} &
k^{\sigma_{i}}c_{ij}\neq1\\
&  & \\
(\xi,\ldots,1,\ldots,\xi)\in\Gamma_{k(i,j),1} & \text{if} & k^{\sigma_{i}%
}c_{ij}=1
\end{array}
\right.  ;
\]
and similarly for the various other cases.

The point is that this is a definable operation on the definable subset
$\Gamma$ of $H^{m}$.

\emph{Bi-interpretation. }The composed monomorphism $H\rightarrow
G\rightarrow\Gamma\subseteq H^{m}$ sends $h\in H$ to%
\begin{align*}
&  (h,1,\ldots,1)~~\ (h\neq1)\\
&  (1,\xi,\ldots,\xi)~~\ (h=1).
\end{align*}
This is definable.

The composed monomorphism $G\rightarrow\Gamma\subseteq H^{m}\subseteq G^{m}$
sends $g\in G$ to $\gamma(g)$. Now%
\[
g\in(H\smallsetminus\{1\})t_{i}\Longleftrightarrow G\models\kappa(gt_{i}%
^{-1})\wedge g\neq t_{i}.
\]
As the map $\gamma$ is definable on each of the sets $(H\smallsetminus
\{1\})t_{i},$ and on the elements $t_{i}$, it follows that $\gamma
:G\rightarrow G^{m}$ is definable.

This completes the proof of \textbf{(b)} when $m\geq3$. \ If $m=2,$ the above
definition gives $\gamma(\xi)=(\xi,1)=\gamma(t_{2}),$ so we need to adjust it.
One way is to replace $\gamma$ with $\gamma^{\ast}:G\rightarrow\Gamma^{\ast
}\subset\Gamma\times\{1,\xi\}\subseteq H^{3},$ where%
\begin{align*}
\gamma^{\ast}(h)  &  =(h,1,1),~\ \gamma^{\ast}(ht_{2})=(1,h,1)\text{ \ \ if
}1\neq h\in H\\
\gamma^{\ast}(t_{2})  &  =(\xi,1,\xi),~~~\gamma^{\ast}(1)=(1,\xi,\xi).
\end{align*}
The group operation is defined in a similar way to the previous case, and
again \textbf{(b)} follows.

Now assume \textbf{(b)}$.$ To each formula $\phi(\mathbf{g};\mathbf{x})$ with
parameters $\mathbf{g}$ from $G$ is associated a formula $\psi(\mathbf{h}%
;\mathbf{y})$ with parameters $\mathbf{h}$ from $H$ and a definable function
$\beta:H^{r}\rightarrow H^{s}$, for some $s$, such that for each
$\mathbf{a}\in H^{r},$%
\[
G\models\phi(\mathbf{g};\mathbf{a})\Longleftrightarrow H\models\psi
(\mathbf{h};\beta(\mathbf{a}))
\]
(cf. \cite{H}, Thm. 5.3.2). Now fix $t\in G$, set $\phi(t;x,y):=(xt=ty)$ and
let $\psi(\mathbf{h};\mathbf{y})$ be the corresponding formula. Then for
$a,b\in H$ we have%
\[
b=t^{-1}at\Longleftrightarrow G\models\phi(t;a,b)\Longleftrightarrow
H\models\psi(\mathbf{h};\beta(a,b)).
\]
Thus the map $a\longmapsto t^{-1}at:H\rightarrow H$ is definable, and
so\textbf{ (a)} holds.
\end{proof}

\bigskip

In many situations, the definability condition \textbf{(a)} is either not
satisfied, or not easy to verify. Under weaker hypotheses, Theorem 3.1 of
\cite{NST} shows that if $H$ is finitely axiomatizable, then so is $G$. The
proof of that theorem, when applied to finitely generated abstract groups,
yields the following simple statement:

\begin{theorem}
Let $\mathcal{C}$ be a class of finitely generated groups, closed under
passing to normal subgroups of finite index. Let $G=\left\langle g_{1}%
,\ldots,g_{d}\right\rangle \in\mathcal{C}$ and let $H=\left\langle
h_{1},\ldots,h_{s}\right\rangle $ be a definable normal subgroup of finite
index in $G$. If $(H,\mathbf{h})$ is finitely axiomatizable in $\mathcal{C}$
then so is $(G,\mathbf{g})$.
\end{theorem}

\noindent (The hypothesis means: $H$ satisfies a sentence $\sigma_{H}(\mathbf{h})$ such
that for any group $\widetilde{H}$ and tuple $\widetilde{\mathbf{h}}%
\in\widetilde{H}^{s}$, if $\widetilde{H}$ satisfies $\sigma_{H}%
(\widetilde{\mathbf{h}})$ then there is an isomorphism $H\rightarrow
\widetilde{H}$ sending $h_{i}$ to $\widetilde{h_{i}}$ for each $i$.)

In view of Theorem \ref{T1}, the converse also holds provided we add the
hypothesis \textbf{(a)}.

With thanks to Andre Nies and Katrin Tent for helpful comments.


\begin{thebibliography}{999}                                                                                              %


\bibitem[H]{H}W. Hodges, \emph{Model Theory}, Cambridge University Press, 1993.

\bibitem[NST]{NST}A. Nies, D. Segal and K. Tent, Finite axiomatizability for
profinite groups, \emph{Proc. London Math. Soc.(3)} \textbf{123} (2021), 597-635.
\end{thebibliography}
\end{document}